# AN OPTIMAL HARDY-MORREY INEQUALITY[*]

GEORGIOS PSARADAKIS

psaradakis.georgios@gmail.com

**Abstract**

In this work we improve the sharp Hardy inequality in the case $p > n$ by adding an optimal weighted Hölder seminorm. To achieve this we first obtain a local improvement. We also obtain a refinement of both the Sobolev inequality for $p > n$ and the Hardy inequality, the latter having the best constant.

# Contents

---

[*]Revised version: contains an extra section with another proof of the alleged optimality

# 1 Introduction

The classic multidimensional version of Hardy's inequality asserts that if $n \geq 1$ is an integer and $p \geq 1$; $p \neq n$, then for all $u \in C_c^\infty(\mathbb{R}^n \setminus \{0\})$

$$\int_{\mathbb{R}^n} |\nabla u|^p dx \geq \Big|\frac{n-p}{p}\Big|^p \int_{\mathbb{R}^n} \frac{|u|^p}{|x|^p} dx, \tag{1.1}$$

where the constant factor $|(n-p)/p|^p$ is sharp; see for example [Mz]-§1.3.1. On the other hand, the classical Sobolev inequality asserts that if $1 < p < n$; $n \geq 2$, then there exists a positive constant $S_{n,p}$, depending only on $n, p$, such that for all $u \in C_c^\infty(\mathbb{R}^n)$

$$\Big(\int_{\mathbb{R}^n} |\nabla u|^p dx\Big)^{1/p} \geq S_{n,p} \Big(\int_{\mathbb{R}^n} |u|^{np/(n-p)} dx\Big)^{(n-p)/np}. \tag{1.2}$$

This inequality is optimal in the sense that it fails if $np/(n-p)$ is replaced by any number $q > np/(n-p)$; see [AF], [GT], [B] and [Mz]. The best constant and the family of extremal functions have been found simultaneously in [A] and [Tl1].

In their pioneering work [BV], Brezis and Vazquez improved Hardy's inequality in the case $p = 2$. They showed that if $\Omega$ is a bounded domain in $\mathbb{R}^n$; $n \geq 3$, containing the origin and $1 < q < 2n/(n-2)$, then there exists a positive constant $C = C(n, q)$ such that for all $u \in C_c^\infty(\Omega)$

$$\Big[\int_\Omega |\nabla u|^2 dx - \Big(\frac{n-2}{2}\Big)^2 \int_\Omega \frac{|u|^2}{|x|^2} dx\Big]^{1/2} \geq \frac{C}{|\Omega|^{1/q-(n-2)/2n}} \Big(\int_\Omega |u|^q dx\Big)^{1/q}. \tag{1.3}$$

In Problem 2 of [BV], the question of whether there is a further improvement in the direction of the inequality (1.3) is posed. An optimal answer was given in [FT], where it was shown that the critical exponent $q = 2n/(n-2)$ is possible after considering a logarithmic correction weight for which the sharp exponent was given. More precisely, it is proved that if $\Omega$ is a bounded domain in $\mathbb{R}^n$; $n \geq 3$, containing the origin, then there exists a positive constant $C = C(n)$ such that for all $u \in C_c^\infty(\Omega)$

$$\Big[\int_\Omega |\nabla u|^2 dx - \Big(\frac{n-2}{2}\Big)^2 \int_\Omega \frac{|u|^2}{|x|^2} dx\Big]^{1/2}$$
$$\geq C\Big(\int_\Omega |u|^{2n/(n-2)} X^{1+n/(n-2)}(|x|/D) dx\Big)^{(n-2)/2n}, \tag{1.4}$$

where $D = \sup_{x \in \Omega} |x|$ and $X(t) = (1 - \log t)^{-1}$, $t \in (0, 1]$. Moreover, the weight function $X^{1+n/(n-2)}$ is optimal in the sense that the power $1 + n/(n-2)$ cannot be decreased; see [AFT] for a second proof where in addition the best constant $C$ is obtained.

Motivated by the above mentioned question of [BV], it is natural to consider the case where $p \neq 2$. Some results in the direction of extending (1.3) and (1.4) in the range $1 < p < n$, were obtained in [ACR]-Theorem 1.1, [BFT]-Theorems B, C and 6.4 and also in [ACP]-Theorem 1.1. In the present paper we focus in the case where $p > n$.

Let us first state Sobolev's inequality in this case: if $p > n \geq 1$ and $\Omega$ is an open set in $\mathbb{R}^n$ with finite volume $|\Omega|$, then there exists a positive constant $s_{n,p}$ depending only on $n, p$, such that

$$\sup_{x \in \Omega} |u(x)| \leq s_{n,p} |\Omega|^{1/n - 1/p} \Big( \int_{\Omega} |\nabla u|^p dx \Big)^{1/p}, \tag{1.5}$$

for all $u \in C_c^{\infty}(\Omega)$. This inequality is also optimal in the sense that it fails if $1/n - 1/p$ is replaced by any exponent $q < 1/n - 1/p$ on $|\Omega|$. The best constant and the family of extremal functions have been found in [Tl2]. C. B. Morrey sharpened this result by replacing the supremum norm with an optimal Hölder semi-norm. What he showed is that there exists a positive constant $C = C(n, p)$ such that

$$\sup_{\substack{x,y \in \mathbb{R}^n \\ x \neq y}} \Big\{ \frac{|u(x) - u(y)|}{|x - y|^{1 - n/p}} \Big\} \leq C \Big( \int_{\mathbb{R}^n} |\nabla u|^p dx \Big)^{1/p}, \tag{1.6}$$

for all $u \in C_c^{\infty}(\mathbb{R}^n)$, and the modulus of continuity $1 - n/p$ is optimal; see [Mr], [AF], [B] and [Mz].

Our main purpose in this paper is to provide optimal Sobolev-type improvements to the sharp Hardy inequality (1.1) for $p > n$. Firstly, we improve both (1.1) and (1.5) by replacing the $L^p$-norm of the length of the gradient in (1.5) with the sharp $L^p$ Hardy difference:

**Theorem A** *Suppose $\Omega$ is a domain in $\mathbb{R}^n$; $n \geq 1$, containing the origin and having finite volume $|\Omega|$. Letting $p > n$, there exists a positive constant $C = C(n, p)$ such that for all $u \in C_c^{\infty}(\Omega \setminus \{0\})$*

$$\sup_{x \in \Omega} |u(x)| \leq C |\Omega|^{1/n - 1/p} \Big[ \int_{\Omega} |\nabla u|^p dx - \Big( \frac{p - n}{p} \Big)^p \int_{\Omega} \frac{|u|^p}{|x|^p} dx \Big]^{1/p}. \tag{1.7}$$

The corresponding result for Morrey's inequality (1.6) is possible only after considering a logarithmic correction weight for which we obtain the sharp exponent. The central result of the paper is the following optimal Hardy-Morrey inequality

**Theorem B** *Suppose $\Omega$ is a bounded domain in $\mathbb{R}^n$; $n \geq 1$, containing the origin and let $p > n$. There exist constants $B = B(n, p) \geq 1$ and $C = C(n, p) > 0$ such that for all $u \in C_c^{\infty}(\Omega \setminus \{0\})$*

$$\sup_{\substack{x,y \in \Omega \\ x \neq y}} \Big\{ \frac{|u(x) - u(y)|}{|x - y|^{1 - n/p}} X^{1/p} \Big( \frac{|x - y|}{D} \Big) \Big\} \leq C \Big[ \int_{\Omega} |\nabla u|^p dx - \Big( \frac{p - n}{p} \Big)^p \int_{\Omega} \frac{|u|^p}{|x|^p} dx \Big]^{1/p}, \tag{1.8}$$

*where $D = B \operatorname{diam}(\Omega)$ and $X(t) = (1 - \log t)^{-1}$; $t \in (0, 1]$. Moreover, the weight function $X^{1/p}$ is optimal in the sense that the power $1/p$ cannot be decreased.*

Note that since $p > n$, one is forced to consider functions in $C_c^{\infty}(\mathbb{R}^n \setminus \{0\})$, i.e. supported away from the origin. This excludes symmetrization techniques as a method of proof. Thus we turn to multidimensional arguments and in particular in Sobolev's integral representation formula. The first step is to show that (1.8) is equivalent to it's counterpart inequality with one point in the Hölder semi-norm taken to be the origin; see Proposition 5.1. In establishing Proposition 5.1, a crucial step is obtaining estimates on balls $B_r$ intersecting $\Omega$ and with arbitrarily small radius. To this end, the following local improvement of the sharp Hardy inequality which is of independent interest, is proved:

**Theorem C** *Suppose $\Omega$ is a bounded domain in $\mathbb{R}^n$; $n \geq 2$, containing the origin and let $p > n$ and $1 \leq q < p$. There exist constants $\Theta = \Theta(n, p, q) \geq 0$ and $C = C(n, p, q) > 0$ such that for all $u \in C_c^\infty(\Omega \setminus \{0\})$, any open ball $B_r$ with $r \in (0, \operatorname{diam}(\Omega))$, and any $D \geq e^\Theta \operatorname{diam}(\Omega)$*

$$r^{n/p} X^{1/p}(r/D) \Big(\frac{1}{|B_r|} \int_{B_r} \frac{|u|^q}{|x|^q} dx\Big)^{1/q} \leq C \Big[\int_\Omega |\nabla u|^p dx - \Big(\frac{p-n}{p}\Big)^p \int_\Omega \frac{|u|^p}{|x|^p} dx\Big]^{1/p}, \qquad (1.9)$$

*where $X(t) = (1 - \log t)^{-1}$; $t \in (0, 1]$.*

The exponent $1/p$ on the logarithmic factor $X^{1/p}$, is translated as the optimal exponent in (1.8). To obtain this exponent in (1.9) we carefully estimate a trace term on the boundary of $B_r$. Let us note here that if one restricts to the family of open balls $B_r$ containing the origin, then Theorem C remains valid for $p < n$ (with the factor $((n-p)/p)^p$ instead of $((p-n)/p)^p$, in the Hardy difference).

We finally note that the second important ingredient in the proof of Theorem B, is to show the optimality of the exponent $1/p$. This is done by finding a suitable family of functions that plays the role of a minimizing sequence for inequality (1.8).

For other directions in strengthening Sobolev's imbedding theorem we refer to [LYZ], [CLYZ] where the $L^p$-norm of the length of the gradient of $u$ in (1.2), (1.5) respectively, was replaced by a smaller quantity, called *the affine energy of $u$*, which is additionally invariant under affine transformations of $\mathbb{R}^n$. We refer also to [CFMP], [C] where the authors improved (1.2), (1.5) respectively, by adding the *functional asymmetry of $u$*, a quantity that measures the distance between $u$ and the family of extremal functions of the inequality (see also [BE]). In [CF], the functional asymmetry of $u$ against the virtual minimizers of (1.1), i.e. the family of functions from which one extracts sequences which imply the sharpness of the best constant $((n-p)/p)^p$, has been added to the right-hand side of (1.1) for any $1 < p < n$. Let us add that in this particular range much progress has been made in case of Hardy's inequality involving distance from the boundary of a given domain; see [Mz], [TT], [BFL], [FMT1], [FMT2] and [FL].

The paper is organized as follows. In §2 we recall some known facts and prove some auxiliary lemmas. The proof of Theorem A and Theorem C is given in §3 and §4 respectively. The multidimensional case of Theorem B is proved in §5, while the case $n = 1$ is treated separately in §6. The optimality of the exponent $1/p$ on $X$ in Theorem B is proved in §7.

# 2 Notation and preparative results

Throughout the paper, $\Omega$ will stand for a domain (i.e. open and connected set) in $\mathbb{R}^n$; $n \geq 1$, containing the origin. Also, we adopt the notation

$$I[u] := \int_\Omega |\nabla u|^p dx - \left|\frac{n-p}{p}\right|^p \int_\Omega \frac{|u|^p}{|x|^p} dx,$$

for $u \in C_c^\infty(\Omega \setminus \{0\})$. By $B_r(x)$ (resp. $\overline{B}_r(x)$) we denote an open (resp. closed) ball in $\mathbb{R}^n$ having radius $r > 0$ and center at $x \in \mathbb{R}^n$. When the center is of no importance we simply write $B_r$. The volume of $B_1$ is denoted by $\omega_n$. We also extend all functions having compact support by zero outside it.

In this section we collect some auxiliary results that will be used throughout the proofs of Theorems A,B and C. In particular we first prove a trace Hardy inequality and then we recall some key lower bounds on $I[u]$ obtained in [BFT]. We conclude with a technical lemma regarding the function $X(t) = (1 - \log t)^{-1}$; $t \in (0, 1]$.

We will need a generalization of (1.1) involving weights and also valid for functions not necessarily vanishing on the boundary of a smooth set $V$.

**Lemma 2.1.** *Let $V$ be a domain in $\mathbb{R}^n$; $n \geq 2$, having locally Lipschitz boundary. Denote by $\nu(x)$ the exterior unit normal vector defined at almost every $x \in \partial V$. For all $q \geq 1$, all $s \neq n$ and any $v \in C_c^\infty(\mathbb{R}^n \setminus \{0\})$, there holds*

$$\int_V \frac{|\nabla v|^q}{|x|^{s-q}} dx + \frac{n-s}{q}\left|\frac{n-s}{q}\right|^{q-2} \int_{\partial V} \frac{|v|^q}{|x|^s} x \cdot \nu(x) dS_x \geq \left|\frac{n-s}{q}\right|^q \int_V \frac{|v|^q}{|x|^s} dx. \qquad (2.1)$$

**Proof.** Integration by parts gives

$$\int_V \nabla |v| \cdot \frac{x}{|x|^s} dx = -\int_V |v| \operatorname{div}\left(\frac{x}{|x|^s}\right) dx + \int_{\partial V} |v| \frac{x}{|x|^s} \cdot \nu dS_x,$$

and since $\operatorname{div}(x|x|^{-s}) = (n-s)|x|^{-s}$ we get

$$\int_V \frac{|\nabla v|}{|x|^{s-1}} dx \geq (s-n) \int_V \frac{|v|}{|x|^s} dx + \int_{\partial V} \frac{|v|}{|x|^s} x \cdot \nu dS_x, \quad \text{if } s > n,$$

$$\int_V \frac{|\nabla v|}{|x|^{s-1}} dx \geq (n-s) \int_V \frac{|v|}{|x|^s} dx - \int_{\partial V} \frac{|v|}{|x|^s} x \cdot \nu dS_x, \quad \text{if } s < n,$$

where we have also used the fact that $|\nabla |v(x)|| \leq |\nabla v(x)|$ for a.e. $x \in V$ (see [LL]-Theorem 6.17). We may write both inequalities in one as follows

$$\int_V \frac{|\nabla v|}{|x|^{s-1}} dx + \frac{n-s}{|n-s|} \int_{\partial V} \frac{|v|}{|x|^s} x \cdot \nu dS_x \geq |n-s| \int_V \frac{|v|}{|x|^s} dx.$$

This is inequality (2.1) for $q = 1$. Substituting $v$ by $|v|^q$ with $q > 1$, we arrive at

$$\frac{q}{|n-s|}\int_V \frac{|\nabla v||v|^{q-1}}{|x|^{s-1}}dx + \frac{n-s}{|n-s|^2}\int_{\partial V}\frac{|v|^q}{|x|^s}x\cdot\nu dS_x \geq \int_V \frac{|v|^q}{|x|^s}dx. \tag{2.2}$$

The first term on the left of (2.2) can be written as follows

$$\begin{aligned}\frac{q}{|n-s|}\int_V \frac{|\nabla v||v|^{q-1}}{|x|^{s-1}}dx &= \int_V \left\{\frac{q}{|n-s|}\frac{|\nabla v|}{|x|^{s/q-1}}\right\}\left\{\frac{|v|^{q-1}}{|x|^{s-s/q}}\right\}dx \\ &\leq \frac{1}{q}\left|\frac{q}{n-s}\right|^q\int_V \frac{|\nabla v|^q}{|x|^{s-q}}dx + \frac{q-1}{q}\int_V \frac{|v|^q}{|x|^s}dx,\end{aligned}$$

by Young's inequality with conjugate exponents $q$ and $q/(q-1)$. Thus (2.2) becomes

$$\frac{1}{q}\left|\frac{q}{n-s}\right|^q\int_V \frac{|\nabla v|^q}{|x|^{s-q}}dx + \frac{n-s}{|n-s|^2}\int_{\partial V}\frac{|v|^q}{|x|^s}x\cdot\nu dS_x \geq \frac{1}{q}\int_V \frac{|v|^q}{|x|^s}dx.$$

Rearranging the constants we arrive at the inequality we sought for. ∎

**Remark 2.2.** A look at the above proof shows that the choice $V = \mathbb{R}^n$ in (2.1) is acceptable provided one cancels the trace term on the left-hand side. More precisely

$$\int_{\mathbb{R}^n} \frac{|\nabla v|^q}{|x|^{s-q}}dx \geq \left|\frac{n-s}{q}\right|^q \int_{\mathbb{R}^n}\frac{|v|^q}{|x|^s}dx, \tag{2.3}$$

for any $q \geq 1$ and all $v \in C_c^\infty(\mathbb{R}^n \setminus \{0\})$. The constant appearing in (2.3) turns out to be sharp, see for example [Mz]-§1.3.1.

Next we quote some known results. In [BFT] the authors obtained various improvements for Hardy's inequality (1.1) valid in any dimension $n \geq 1$. In particular, the substitution $u(x) = |x|^{1-n/p}v(x)$ and elementary vectorial inequalities lead to the following lower bounds on $I[u]$; see [BFT]-Lemma 3.3,

$$I[u] \geq c_p \int_\Omega |x|^{p-n}|\nabla v|^p dx, \tag{2.4}$$

$$I[u] \geq c_p \left|\frac{p-n}{p}\right|^{p-2}\int_\Omega |x|^{2-n}|v|^{p-2}|\nabla v|^2 dx, \tag{2.5}$$

both in case $n \neq p \geq 2$. Here $c_p$ is a positive constant depending only on $p$. They also obtained the optimal homogeneous improvement to Hardy's inequality:

**Theorem 2.3** ([BFT]). *Let $n \neq p > 1$ and suppose $\Omega$ is a bounded domain in $\mathbb{R}^n$ containing the origin. There exists a constant $\theta \geq 0$ depending only on $n, p$, such that for all $u \in C_c^\infty(\Omega\setminus\{0\})$ and any $D \geq e^\theta \sup_{x\in\Omega}|x|$*

$$I[u] \geq \frac{p-1}{2p}\left|\frac{p-n}{p}\right|^{p-2}\int_\Omega \frac{|u(x)|^p}{|x|^p}X^2(|x|/D)dx,$$

*where $X(t) = (1-\log t)^{-1}$, $t \in (0,1]$. The weight function $X^2$ is optimal, in the sense that the power 2 cannot be decreased, and the constant appearing on the right-hand side is sharp.*

Let us note that the one dimensional case with $p = 2$ of the above theorem was treated in the appendix of [BM].

The next lemma is the special case $n = 1$ of [BFT]-Proposition 3.4.

**Lemma 2.4.** *Let $1 < p < 2$. There exists a positive constant $c = c(p)$ such that for all $u \in C_c^\infty(0, R)$ and any $D \geq R$, the following inequality is valid with $v = t^{1/p-1}u$*

$$I[u] \geq c \int_0^R t^{p-1} |v'(t)|^p X^{2-p}(t/D) dt. \tag{2.6}$$

We close this section with the following technical fact concerning the auxiliary function $X(t) = (1 - \log t)^{-1}$, $t \in (0, 1]$, which will be helpful in our computations.

**Lemma 2.5.** *Let $\alpha > -1$ and $\beta, R > 0$. For all $r \in (0, R]$, all $c > 1/(\alpha+1)$ and any $D \geq e^\eta R$, where $\eta := \max\{0, \frac{(\beta-\alpha-1)c+1}{(\alpha+1)c-1}\}$, we have*

$$(i) \quad \int_0^r t^\alpha X^{-\beta}(t/D) dt \leq c r^{\alpha+1} X^{-\beta}(r/D).$$

*If $\alpha$ is restricted in $(-1, 0]$ then for all $0 \leq y \leq x \leq R$ we have*

$$(ii) \quad \int_y^x t^\alpha X^{-\beta}(t/D) dt \leq c(x-y)^{\alpha+1} X^{-\beta}((x-y)/D).$$

**Proof.** Let $c > 0$ and $D \geq R$. We set

$$f(r) := \int_0^r t^\alpha X^{-\beta}(t/D) dt - c r^{\alpha+1} X^{-\beta}(r/D), \quad r \in (0, R].$$

To prove $(i)$ it suffices to show that for suitable values of the parameters $c$ and $D$, we have $f(r) \leq 0$ for all $r \in (0, R)$. We have $f(0+) = 0$ and thus it is enough to choose $c$ and $D$ in such a way so that $f$ is decreasing in $(0, R)$. To this end we compute

$$f'(r) = c r^\alpha X^{-\beta}(r/D)[1/c - (\alpha+1) + \beta X(r/D)], \quad r \in (0, R].$$

It is easy to see that for $c > 1/(\alpha+1)$, any $D \geq e^\eta R$ is such that $f'(r) \leq 0$ for all $r \in (0, R)$. To prove $(ii)$ we note that since $f$ is decreasing, $0 \leq y \leq x \leq R$ implies $f(y) \geq f(x)$, and so

$$\int_y^x t^\alpha X^{-\beta}(t/D) dt \leq c[x^{\alpha+1} X^{-\beta}(x/D) - y^{\alpha+1} X^{-\beta}(y/D)]$$
$$\leq c(x^{\alpha+1} - y^{\alpha+1}) X^{-\beta}(x/D)$$
$$\leq c(x^{\alpha+1} - y^{\alpha+1}) X^{-\beta}((x-y)/D),$$

where the last two inequalities follow since $X^{-\beta}(r/D)$ is decreasing in $(0, R)$. If $\alpha \in (-1, 0]$ then $x^{\alpha+1} - y^{\alpha+1} \leq (x-y)^{\alpha+1}$, and the result follows. ∎

# 3 The Hardy-Sobolev inequality for $p > n \geq 1$

In this section we prove Theorem A. Assume first that $n \geq 2$. Letting $x \in \Omega$, by the standard representation formula (see [GT]-Lemma 7.14) we have

$$
\begin{aligned}
u(x) &= \frac{1}{n\omega_n} \int_\Omega \frac{(x-z) \cdot \nabla u(z)}{|x-z|^n} dz \\
&\leq \frac{1}{n\omega_n} \int_\Omega \frac{|\nabla u(z)|}{|x-z|^{n-1}} dz.
\end{aligned} \tag{3.1}
$$

Setting $u(z) = |z|^{1-n/p} v(z)$, we arrive at

$$
|u(x)| \leq \frac{1}{n\omega_n} \underbrace{\int_\Omega \frac{|z|^{1-n/p}|\nabla v(z)|}{|x-z|^{n-1}} dz}_{=:K(x)} + \frac{p-n}{pn\omega_n} \underbrace{\int_\Omega \frac{|v(z)|}{|z|^{n/p}|x-z|^{n-1}} dz}_{=:\Lambda(x)}. \tag{3.2}
$$

We first estimate $K(x)$. Using Hölder's inequality we get

$$
K(x) \leq \Big( \underbrace{\int_\Omega \frac{1}{|x-z|^{(n-1)p/(p-1)}} dz}_{=:M(x)} \Big)^{1-1/p} \Big( \int_\Omega |z|^{p-n} |\nabla v(z)|^p dz \Big)^{1/p}. \tag{3.3}
$$

Note that $M(x)$ is finite since $(n-1)p/(p-1) < n$ if and only if $p > n$. To estimate $M(x)$ we set $R := (|\Omega|/\omega_n)^{1/n}$, so that the volume of a ball with radius $R$ to be equal to the volume of $\Omega$. Then $M(x)$ increases if we change the domain of integration from $\Omega$ to $B_R(x)$. Therefore

$$
\begin{aligned}
M(x) &\leq \int_{B_R(x)} \frac{1}{|x-z|^{(n-1)p/(p-1)}} dz \\
&= n\omega_n \frac{p-1}{p-n} (|\Omega|/\omega_n)^{(p-n)/n(p-1)},
\end{aligned}
$$

and using (2.4) in the second factor of (3.3), we get

$$
K(x) \leq C_1(n,p) |\Omega|^{1/n - 1/p} (I[u])^{1/p}. \tag{3.4}
$$

Next we bound $\Lambda(x)$. Using Hölder's inequality with conjugate exponents $p/(p-1-\varepsilon)$ and $p/(1+\varepsilon)$, where $0 < \varepsilon < (p-n)/n$ is fixed and depending only on $n, p$, we get

$$
\Lambda(x) \leq \Big( \underbrace{\int_\Omega \frac{1}{|x-z|^{(n-1)p/(p-1-\varepsilon)}} dz}_{=:M(\varepsilon,x)} \Big)^{1-(1+\varepsilon)/p} \Big( \int_\Omega \frac{|v(z)|^{p/(1+\varepsilon)}}{|z|^{n/(1+\varepsilon)}} dz \Big)^{(1+\varepsilon)/p}. \tag{3.5}
$$

Note that $M(\varepsilon, x)$ is finite since $(n-1)p/(p-1-\varepsilon) < n$ if and only if $\varepsilon < (p-n)/n$. As before, we set $R := (|\Omega|/\omega_n)^{1/n}$ and $M(\varepsilon, x)$ increases if we change the domain of integration from $\Omega$ to $B_R(x)$. Therefore,

$$
\begin{aligned}
M(\varepsilon, x) &\leq \int_{B_R(x)} \frac{1}{|x-z|^{(n-1)p/(p-1-\varepsilon)}} dz \\
&= n\omega_n \frac{p-1-\varepsilon}{p-n-n\varepsilon}(|\Omega|/\omega_n)^{(p-n-n\varepsilon)/n(p-1-\varepsilon)},
\end{aligned}
$$

and using (2.3) for $s = n/(1+\varepsilon)$ and $q = p/(1+\varepsilon)$ in the second factor of (3.5), we obtain

$$
\Lambda(x) \leq C_2(n,p)|\Omega|^{1/n-1/p-\varepsilon/p}\Big(\int_\Omega |z|^{(p-n)/(1+\varepsilon)}|\nabla v(z)|^{p/(1+\varepsilon)} dz\Big)^{(1+\varepsilon)/p}.
$$

Using once more Hölder's inequality with conjugate exponents $1 + 1/\varepsilon$ and $1 + \varepsilon$, we get

$$
\begin{aligned}
\Lambda(x) &\leq C_2(n,p)|\Omega|^{1/n-1/p-\varepsilon/p}\Big[|\Omega|^{\varepsilon/(1+\varepsilon)}\Big(\int_\Omega |z|^{p-n}|\nabla v(z)|^p dz\Big)^{1/(1+\varepsilon)}\Big]^{(1+\varepsilon)/p} \\
&= C_2(n,p)|\Omega|^{1/n-1/p}\Big(\int_\Omega |z|^{p-n}|\nabla v(z)|^p dz\Big)^{1/p}. \\
\text{(by (2.4))} &\leq C_3(n,p)|\Omega|^{1/n-1/p}(I[u])^{1/p}. \quad (3.6)
\end{aligned}
$$

The proof follows inserting (3.6) and (3.4) in (3.2). For $n = 1$ we have

$$
\begin{aligned}
u(x) &\leq \frac{1}{2}\int_0^R |u'(t)| dt \\
(\text{setting } u(t) = t^{1-1/p} v(t)) &\leq \frac{1}{2}\int_0^R t^{1-1/p}|v'(t)| dt + \frac{p-1}{2p}\int_0^R t^{-1/p}|v(t)| dt \\
&\leq \int_0^R t^{1-1/p}|v'(t)| dt,
\end{aligned}
$$

by (2.3) for $n = q = 1$ and $s = 1/p$. The proof follows applying Hölder's inequality and using (2.4) for $p \geq 2$ or (2.6) if $1 < p < 2$. We omit further details. ∎

## 4  A local improvement

Here we prove Theorem C in case $0 \in \overline{B}_r$. We emphasize that under this assumption we will prove (1.9) for general $p, q > 1$ satisfying $1 \leq q < p$, $p \neq n$. The proof in case $0 \notin \overline{B}_r$ and $p > n$ is given after the proof of Proposition 5.1 of the next section.

**Proof of Theorem C in case $0 \in \overline{B}_r$.** Let $\Omega$ be a bounded domain in $\mathbb{R}^n$ containing the origin and let $1 \leq q < p$, $p \neq n$. Suppose $u \in C_c^\infty(\Omega \setminus \{0\})$ and let also $B_r$ be any ball containing zero with $r \in (0, \text{diam}(\Omega))$. Setting $u(x) = |x|^{1-n/p}v(x)$ we get

$$\int_{B_r} \frac{|u|^q}{|x|^q}dx = \int_{B_r} \frac{|v|^q}{|x|^{nq/p}}dx$$

$$\leq \left[\frac{pq}{n(p-q)}\right]^q \underbrace{\int_{B_r} |x|^{q(p-n)/p}|\nabla v|^q dx}_{=:N_r}$$

$$+ \frac{pq}{n(p-q)} \underbrace{\int_{\partial B_r} \frac{|v|^q}{|x|^{nq/p}} x \cdot \nu dS_x}_{=:P_r}, \quad (4.1)$$

where we have used Lemma 2.1 for $V = B_r$ and $s = nq/p$. We use Hölder's inequality with conjugate exponents $p/(p-q)$ and $p/q$, to get

$$N_r \leq (\omega_n r^n)^{(p-q)/p} \left(\int_{B_r} |x|^{p-n}|\nabla v|^p dx\right)^{q/p}$$

$$\text{(by (2.4))} \leq C_1(n,p,q) r^{n(p-q)/p}(I[u])^{q/p}$$

$$\leq C_1(n,p,q) r^{n(p-q)/p} X^{-q/p}(r/D)(I[u])^{q/p}, \quad (4.2)$$

for any $D \geq \text{diam}(\Omega)$ since $0 \leq X(t) \leq 1$ for all $t \in (0,1]$. For $P_r$ we write first

$$P_r = \int_{\partial B_r} \left\{X^{-q/p}(|x|/D)\right\}\left\{\frac{|v|^q}{|x|^{nq/p}}X^{q/p}(|x|/D)\right\} x \cdot \nu dS_x$$

$$\leq \underbrace{\left(\int_{\partial B_r} X^{-q/(p-q)}(|x|/D) x \cdot \nu dS_x\right)^{(p-q)/p}}_{=:S_r} \underbrace{\left(\int_{\partial B_r} \frac{|v|^p}{|x|^n} X(|x|/D) x \cdot \nu dS_x\right)^{q/p}}_{=:T_r}, (4.3)$$

where we have used once more Hölder's inequality with exponents $p/(p-q)$ and $p/q$. Note that $\nu$ is the outward pointing unit normal vector field along $\partial B_r$ and that since $0 \in \overline{B}_r$ we have $x \cdot \nu \geq 0$ for all $x \in \partial B_r$. By the divergence theorem we have

$$S_r = \int_{B_r} \text{div}[X^{-q/(p-q)}(|x|/D) \, x]dx$$

$$= n \int_{B_r} X^{-q/(p-q)}(|x|/D)dx - \frac{q}{p-q}\int_{B_r} X^{1-q/(p-q)}(|x|/D)dx$$

$$\leq n \int_{B_r(0)} X^{-q/(p-q)}(|x|/D)dx,$$

since the integral increases if we change the domain of integration from $B_r$ to $B_r(0)$. Thus

$$S_r \leq n^2 \omega_n \int_0^r t^{n-1} X^{-q/(p-q)}(t/D)dt$$

$$\leq C_2(n) r^n X^{-q/(p-q)}(r/D), \quad (4.4)$$

for any $D \geq e^\eta \operatorname{diam}(\Omega)$, where $\eta \geq 0$ depends only on $n, p, q$, due to Lemma 2.5 for $\alpha = n-1$ and $\beta = q/(p-q)$. $T_r$ will be estimated after an integration by parts. More precisely

$$T_r = \int_{B_r} \operatorname{div}\left\{\frac{X(|x|/D)}{|x|^n}x\right\}|v|^p dx + \int_{B_r} \frac{X(|x|/D)}{|x|^n}x \cdot \nabla(|v|^p)dx.$$

A simple calculation shows that $\operatorname{div}\{\frac{X(|x|/D)}{|x|^n}x\} = \frac{X^2(|x|/D)}{|x|^n}$ for any $x \in \Omega \setminus \{0\}$. In the second integral we compute the gradient and note that $x \cdot \nabla |v(x)| \leq |x||\nabla v(x)|$ for a.e. $x \in \Omega$. Thus,

$$T_r \leq \int_\Omega \frac{|v|^p}{|x|^n}X^2(|x|/D)dx + p\int_\Omega |x|^{1-n}|v|^{p-1}|\nabla v|X(|x|/D)dx.$$

We rearrange the integrand in the second integral above as follows

$$\begin{aligned} T_r &\leq \int_\Omega \frac{|v|^p}{|x|^n}X^2(|x|/D)dx + p\int_\Omega \left\{|x|^{1-n/2}|v|^{p/2-1}|\nabla v|\right\}\left\{\frac{|v|^{p/2}}{|x|^{n/2}}X(|x|/D)\right\}dx \\ &\leq \int_\Omega \frac{|v|^p}{|x|^n}X^2(|x|/D)dx + p\left(\int_\Omega |x|^{2-n}|v|^{p-2}|\nabla v|^2 dx\right)^{1/2}\left(\int_\Omega \frac{|v|^p}{|x|^n}X^2(|x|/D)dx\right)^{1/2}, \end{aligned}$$

by the Cauchy-Schwarz inequality. According to Theorem 2.3, there exist constants $\theta \geq 0$ and $b > 0$, both depending only on $n, p$, such that for any $D \geq e^\theta \sup_{x \in \Omega}|x|$, the first term and the second radicand on the right hand side (when returned to the original function by $v(x) = |x|^{n/p-1}u(x)$) are bounded above by $bI[u]$. Due to (2.5), the first radicand is also bounded above by $C(n,p)I[u]$. It follows that

$$T_r \leq C_3(n,p)I[u], \tag{4.5}$$

for any $D \geq e^\theta \sup_{x \in \Omega}|x|$. Setting $\Theta = \max\{\theta, \eta\}$ and noting that $0 \in \Omega$ implies $\sup_{x \in \Omega}|x| \leq \operatorname{diam}(\Omega)$, we insert (4.4) and (4.5) into estimate (4.3) to obtain

$$P_r \leq C_4(n,p,q)r^{n(p-q)/p}X^{-q/p}(r/D)(I[u])^{q/p},$$

for any $D \geq e^\Theta \operatorname{diam}(\Omega)$. The last inequality together with (4.2), when applied to estimate (4.1), give

$$\int_{B_r} \frac{|u|^q}{|x|^q}dx \leq C_5(n,p,q)r^{n(p-q)/p}X^{-q/p}(r/D)(I[u])^{q/p}.$$

for any $D \geq e^\Theta \operatorname{diam}(\Omega)$. Rearranging, raising in $1/q$ and taking the supremum over all $B_r$ containing zero with $r \in (0, \operatorname{diam}(\Omega))$, the result follows. ∎

# 5 The multidimensional Hardy-Morrey inequality

In this section we prove Theorem B when $n \geq 2$. We first obtain (1.8) with one point in the Hölder semi-norm taken to be the origin. More precisely, we prove

**Proposition 5.1.** Let $p > n \geq 2$ and suppose $\Omega$ is a bounded domain in $\mathbb{R}^n$ containing the origin. There exist constants $\Theta \geq 0$ and $C > 0$ both depending only on $n, p$, such that for all $u \in C_c^\infty(\Omega \setminus \{0\})$ and any $D \geq e^\Theta \operatorname{diam}(\Omega)$

$$\sup_{x \in \Omega} \left\{ \frac{|u(x)|}{|x|^{1-n/p}} X^{1/p}(|x|/D) \right\} \leq C(I[u])^{1/p}. \tag{5.1}$$

**Proof.** Let $B_r$ be a ball containing zero with $r \in (0, \operatorname{diam}(\Omega))$ and set $u_{B_r} = |B_r|^{-1} \int_{B_r} u \, dz$. Letting $x \in B_r$, the local version of the representation formula (3.1) ([GT]-Lemma 7.16), asserts

$$|u(x) - u_{B_r}| \leq \frac{2^n}{n \omega_n} \int_{B_r} \frac{|\nabla u(z)|}{|x-z|^{n-1}} dz.$$

Setting $u(z) = |z|^{1-n/p} v(z)$, we arrive at

$$\frac{n\omega_n}{2^n} |u(x) - u_{B_r}| \leq \underbrace{\int_{B_r} \frac{|z|^{1-n/p} |\nabla v(z)|}{|x-z|^{n-1}} dz}_{=: K_r(x)} + \frac{p-n}{p} \underbrace{\int_{B_r} \frac{|v(z)|}{|z|^{n/p} |x-z|^{n-1}} dz}_{=: \Lambda_r(x)}. \tag{5.2}$$

We will derive suitable bounds for $K_r(x), \Lambda_r(x)$. For $K_r(x)$ we use Hölder's inequality

$$K_r(x) \leq \left( \int_{B_r} \frac{1}{|x-z|^{(n-1)p/(p-1)}} dz \right)^{1-1/p} \left( \int_{B_r} |z|^{p-n} |\nabla v|^p dz \right)^{1/p}.$$

Both integrals increase if we integrate over $B_r(x)$ and $\Omega$ respectively. Hence

$$K_r(x) \leq \left( \int_{B_r(x)} \frac{1}{|x-z|^{(n-1)p/(p-1)}} dz \right)^{1-1/p} \left( \int_\Omega |z|^{p-n} |\nabla v|^p dz \right)^{1/p}.$$

Computing the first integral and using (2.4) for the second, we arrive at

$$\begin{aligned} K_r(x) &\leq C_1(n,p) r^{1-n/p} (I[u])^{1/p} \\ &\leq C_1(n,p) r^{1-n/p} X^{-1/p}(r/D) (I[u])^{1/p}, \end{aligned} \tag{5.3}$$

for any $D \geq \operatorname{diam}(\Omega)$, where the last inequality follows since $0 < X(t) \leq 1$ for all $t \in (0, 1]$. Next we bound $\Lambda_r(x)$. Using Hölder's inequality with conjugate exponents $p/(p-1-\varepsilon)$ and $p/(1+\varepsilon)$, where $0 < \varepsilon < (p-n)/n$ is fixed but depending only on $n, p$, we obtain

$$\Lambda_r(x) \leq \underbrace{\left( \int_{B_r} \frac{1}{|x-z|^{(n-1)p/(p-1-\varepsilon)}} dz \right)^{1-(1+\varepsilon)/p}}_{=: M_r(x)} \left( \int_{B_r} \frac{|v|^{p/(1+\varepsilon)}}{|z|^{n/(1+\varepsilon)}} dz \right)^{(1+\varepsilon)/p}. \tag{5.4}$$

Note that $M_r(x)$ is finite since $(n-1)p/(p-1-\varepsilon) < n$ if and only if $\varepsilon < (p-n)/n$. Moreover, we note that it increases if we change the domain of integration from $B_r$ to $B_r(x)$. Therefore,

$$\begin{aligned} M_r(x) &\leq \int_{B_r(x)} \frac{1}{|x-z|^{(n-1)p/(p-1-\varepsilon)}} dz \\ &= C_2(n,p) r^{(p-n-n\varepsilon)/(p-1-\varepsilon)}, \end{aligned}$$

and returning to the original function on the second integral on the right of (5.4), we obtain

$$\Lambda_r(x) \leq C_2(n,p)r^{1-n/p-n\varepsilon/p}\left(\int_{B_r}\frac{|u|^{p/(1+\varepsilon)}}{|z|^{p/(1+\varepsilon)}}dz\right)^{(1+\varepsilon)/p}.$$

Using Theorem C with $q = p/(1+\varepsilon)$,

$$\begin{aligned}\Lambda_r(x) &\leq C_3(n,p)r^{1-n/p-n\varepsilon/p}\left(r^{n\varepsilon/(1+\varepsilon)}X^{-1/(1+\varepsilon)}(r/D)(I[u])^{1/(1+\varepsilon)}\right)^{(1+\varepsilon)/p} \\ &= C_3(n,p)r^{1-n/p}X^{-1/p}(r/D)(I[u])^{1/p},\end{aligned} \quad (5.5)$$

for any $D \geq e^{\Theta}\operatorname{diam}(\Omega)$, where $\Theta$ depends only on $n, p, \varepsilon$ (and thus only on $n, p$).

Applying estimates (5.5) and (5.3) to estimate (5.2), we finally conclude

$$|u(x) - u_{B_r}| \leq C_4(n,p)r^{1-n/p}X^{-1/p}(r/D)(I[u])^{1/p}, \quad (5.6)$$

for all $x \in B_r$ and any $D \geq e^{\Theta}\operatorname{diam}(\Omega)$. Since $0 \in B_r$, it follows from (5.6), that

$$|u_{B_r}| \leq C_4(n,p)r^{1-n/p}X^{-1/p}(r/D)(I[u])^{1/p}.$$

Hence

$$\begin{aligned}|u(x)| &\leq |u(x) - u_{B_r}| + |u_{B_r}| \\ &\leq 2C_4(n,p)r^{1-n/p}X^{-1/p}(r/D)(I[u])^{1/p},\end{aligned}$$

for all $x \in B_r$ and any $D \geq e^{\Theta}\operatorname{diam}(\Omega)$. Now if $x \in \Omega$, we consider a ball $B_r$ of radius $r = |x|$, containing $x$. Then the previous inequality yields

$$|u(x)| \leq C(n,p)|x|^{1-n/p}X^{-1/p}(|x|/D)(I[u])^{1/p},$$

for any $D \geq e^{\Theta}\operatorname{diam}(\Omega)$. Rearranging and taking the supremum over all $x \in \Omega$, the result follows. ∎

**Completion of proof of Theorem C.** Let $r \in (0, \operatorname{diam}(\Omega))$ and $p > n$. Using (5.1) we obtain

$$\begin{aligned}\int_{B_r}\frac{|u|^q}{|x|^q}dx &\leq C^q(I[u])^{q/p}\int_{B_r}\frac{1}{|x|^{nq/p}X^{q/p}(|x|/D)}dx \\ &\leq C^q(I[u])^{q/p}\int_{B_r(0)}\frac{1}{|x|^{nq/p}X^{q/p}(|x|/D)}dx \\ &= C^q n\omega_n(I[u])^{q/p}\int_0^r t^{n-1-nq/p}X^{-q/p}(t/D)dt \\ &\leq C(n,q,p)r^{n-nq/p}X^{-q/p}(r/D)(I[u])^{q/p},\end{aligned}$$

for any $D \geq \max\{e^{\eta}, e^{\Theta}\}\operatorname{diam}(\Omega)$, where $\eta = \eta(n,p,q)$, by Lemma 2.5. Rearranging, raising in $1/q$ and taking the supremum over all $B_r$ with $r \in (0, \operatorname{diam}(\Omega))$, we arrive at (1.9) without having assumed $0 \in \overline{B}_r$. ∎

Now we utilize (5.1) in order to obtain its counterpart inequality with the exact Hölder seminorm, i.e. inequality (1.8).

**Proof of Theorem B in case $n \geq 2$.** Letting $x, y \in \Omega$ with $x \neq y$, we consider a ball $B_r$ of radius $r := |x - y|$ containing $x, y$. Note that $r \in (0, \text{diam}(\Omega))$. We have

$$\begin{aligned}
|u(x) - u(y)| &\leq |u(x) - u_{B_r}| + |u(y) - u_{B_r}| \\
&\leq \frac{2^n}{n\omega_n}\Big\{\underbrace{\int_{B_r} \frac{|\nabla u(z)|}{|x - z|^{n-1}}dz}_{=:J(x)} + \underbrace{\int_{B_r} \frac{|\nabla u(z)|}{|y - z|^{n-1}}dz}_{=:J(y)}\Big\},
\end{aligned} \quad (5.7)$$

where we have used [GT]-Lemma 7.16 twice. We will bound $J(x)$ independently on $x$ so that the same estimate holds also for $J(y)$. We start with the substitution $u(z) = |z|^{1-n/p}v(z)$, to get

$$J(x) \leq \underbrace{\int_{B_r} \frac{|z|^{1-n/p}|\nabla v(z)|}{|x - z|^{n-1}}dz}_{=:K_r(x)} + \frac{p-n}{p}\underbrace{\int_{B_r} \frac{|v(z)|}{|z|^{n/p}|x - z|^{n-1}}dz}_{=:\Lambda_r(x)}. \quad (5.8)$$

We estimate $K_r(x)$ in the same manner as in (5.3). So

$$K_r(x) \leq C_1(n,p) r^{1-n/p} X^{-1/p}(r/D)(I[u])^{1/p}, \quad (5.9)$$

for any $D \geq \text{diam}(\Omega)$. To estimate $\Lambda_r(x)$ we return to the original function by $v(z) = |z|^{n/p-1}u(z)$, thus

$$\Lambda_r(x) = \int_{B_r} \frac{|u(z)|}{|z||x - z|^{n-1}}dz.$$

Inserting (5.1) in $\Lambda_r(x)$, we obtain

$$\Lambda_r(x) \leq C_2(n,p)(I[u])^{1/p}\int_{B_r} \frac{X^{-1/p}(|z|/D)}{|z|^{n/p}|x - z|^{n-1}}dz,$$

for any $D \geq e^\Theta \text{diam}(\Omega)$. Using Hölder's inequality with conjugate exponents $Q$ and $Q' = Q/(Q-1)$, where $n < Q < p$ is fixed but depending only on $n, p$, we deduce

$$\Lambda_r(x) \leq C_2(n,p)(I[u])^{1/p}\Big(\int_{B_r} \frac{X^{-Q/p}(|z|/D)}{|z|^{nQ/p}}dz\Big)^{1/Q}\Big(\int_{B_r} \frac{1}{|x - z|^{(n-1)Q'}}dz\Big)^{1/Q'}.$$

Note that both integrals above are finite since $nQ/p < n$ if and only if $Q < p$, and $(n-1)Q' < n$ if and only if $n < Q$. Further, both integrals increase if we integrate over $B_r(0)$ and $B_r(x)$ respectively. Therefore

$$\begin{aligned}
\Lambda_r(x) &\leq C_2(n,p)(I[u])^{1/p}\Big(\int_{B_r(0)} \frac{X^{-Q/p}(|z|/D)}{|z|^{nQ/p}}dz\Big)^{1/Q}\Big(\int_{B_r(x)} \frac{1}{|x - z|^{(n-1)Q'}}dz\Big)^{1/Q'} \\
&= C_3(n,p)(I[u])^{1/p}\Big(\int_0^r t^{n-1-nQ/p}X^{-Q/p}(t/D)dt\Big)^{1/Q} r^{n/Q'-n+1},
\end{aligned} \quad (5.10)$$

for any $D \geq e^\Theta \operatorname{diam}(\Omega)$. Lemma 2.5 for $\alpha = n-1-nQ/p$ and $\beta = Q/p$ ensures the existence of constants $\eta \geq 0$ and $c > 0$ both depending only on $n, p, Q$ (and thus only on $n, p$), such that

$$\int_0^r t^{n-1-nQ/p} X^{-Q/p}(t/D) dt \leq c r^{n-nQ/p} X^{-Q/p}(r/D),$$

for any $D \geq e^\eta \operatorname{diam}(\Omega)$. Thus (5.10) becomes

$$\Lambda_r(x) \leq C_4(n,p) r^{1-n/p} X^{-1/p}(r/D) (I[u])^{1/p}, \tag{5.11}$$

for any $D \geq e^{\Theta'} \operatorname{diam}(\Omega)$ where $\Theta' = \max\{\Theta, \eta\}$.

Summarizing, in view of (5.9) and (5.11), estimate (5.8) becomes

$$J(x) \leq C_5(n,p) r^{1-n/p} X^{-1/p}(r/D) (I[u])^{1/p},$$

for any $D \geq e^{\Theta'} \operatorname{diam}(\Omega)$. The same estimate holds for $J(y)$ and thus (5.7) becomes

$$|u(x) - u(y)| \leq C_6(n,p) r^{1-n/p} X^{-1/p}(r/D) (I[u])^{1/p},$$

for any $D \geq e^{\Theta'} \operatorname{diam}(\Omega)$. The proof of (1.8) in case $n \geq 2$ is completed with $B = e^{\Theta'}$. ∎

## 6 The one-dimensional Hardy-Morrey inequality

Theorem B in the one-dimensional case has an easier proof. We present it in this separate section. Firstly note that it suffices to restrict ourselves in the case $\Omega = (0, R); R > 0$. We start with a lemma which can be obtained by [Mz]-§1.3.2-Theorem 2. We give the proof for convenience.

**Lemma 6.1.** *Let $q > 1$, $\beta > 1 - q$. There exists a positive constant $c = c(q, \beta)$ such that, for all $v \in C_c^\infty(0, R)$ and any $D \geq R$*

$$\sup_{x \in (0,R)} \left\{ |v(x)| X^{(\beta+q-1)/q}(x/D) \right\} \leq c \left( \int_0^R t^{q-1} |v'(t)|^q X^\beta(t/D) dt \right)^{1/q}.$$

**Proof.** Letting $D \geq R$, we have

$$\begin{aligned}
v(x) &= -\int_x^D v'(t) dt \\
&\leq \left( \int_x^D t^{-1} X^{-\beta/(q-1)}(t/D) dt \right)^{1-1/q} \left( \int_x^D t^{q-1} |v'(t)|^q X^\beta(t/D) dt \right)^{1/q},
\end{aligned}$$

where in the last inequality we have used Hölder's inequality with conjugate exponents $q$ and $q/(q-1)$. Since $v \in C_c^\infty(0, R)$ the second integral is actually over $(x, R)$. In addition

$$(X^{-1-\beta/(q-1)}(t/D))' = (-1 - \beta/(q-1)) t^{-1} X^{-\beta/(q-1)}(t/D),$$

so that the left integral can be computed. We find

$$|v(x)| \leq c\Big(X^{-1-\beta/(q-1)}(x/D) - 1\Big)^{1-1/q} \Big(\int_0^R t^{q-1}|v'(t)|^q X^\beta(t/D)dt\Big)^{1/q},$$

where $c = [(q-1)/(q-1+\beta)]^{1-1/q}$, or,

$$|v(x)|X^{(\beta+q-1)/q}(x/D) \leq c\Big(1-X^{(\beta+q-1)/(q-1)}(x/D)\Big)^{1-1/q} \Big(\int_0^R t^{q-1}|v'(t)|^q X^\beta(t/D)dt\Big)^{1/q}.$$

The result follows since $(1 - X^{(\beta+q-1)/(q-1)}(x/D))^{1-1/q} \leq 1$, for all $x \in (0, R]$. ∎

In correspondence to the case $n \geq 2$, we proceed by proving Theorem B with one point in the Hölder semi-norm taken to be the origin. In particular, we have

**Proposition 6.2.** Let $p > 1$. There exists positive constant $c_p$ depending only on $p$, such that for all $u \in C_c^\infty(0, R)$ and any $D \geq R$

$$\sup_{x \in (0,R)} \Big\{\frac{|u(x)|}{x^{1-1/p}} X^{1/p}(x/D)\Big\} \leq c_p (I[u])^{1/p}. \tag{6.1}$$

**Proof.** We set $v(x) = x^{-1+1/p} u(x)$. If $1 < p < 2$, by Lemma 6.1 for $q = p$ and $\beta = 2 - p$, we have

$$|v(x)|X^{1/p}(x/D) \leq c_p \Big(\int_0^R t^{p-1}|v'(t)|^p X^{2-p}(t/D)dt\Big)^{1/p},$$

for any $D \geq R$. The result follows by (2.6). If $p \geq 2$, by Lemma 6.1 for $q = 2$ and $\beta = 0$, we have

$$|w(x)|X^{1/2}(x/D) \leq c_p \Big(\int_0^R t|w'(t)|^2 dt\Big)^{1/2},$$

for any $w \in C_c^\infty(0, R)$ and any $D \geq R$. For $w(x) = |v(x)|^{p/2}$, we obtain

$$|v(x)|X^{1/p}(x/D) \leq c_p \Big(\int_0^R t|v(t)|^{p-2}|v'(t)|^2 dt\Big)^{1/p}.$$

The result follows by (2.5). ∎

Now we use (6.1) to obtain its counterpart inequality with the exact Hölder semi-norm.

**Proof of Theorem B in case $n = 1$.** For $0 < y < x < R$ we have

$$|u(x) - u(y)| = \Big|\int_y^x u'(t)dt\Big|$$

$$(\text{setting } u(t) = t^{1-1/p}v(t)) \leq \underbrace{\int_y^x t^{1-1/p}|v'(t)|dt}_{:=K(x,y)} + \frac{p-1}{p}\underbrace{\int_y^x \frac{|v(t)|}{t^{1/p}}dt}_{:=\Lambda(x,y)}. \tag{6.2}$$

To estimate $\mathrm{K}(x,y)$ we use Hölder's inequality to get

$$\begin{aligned}
\mathrm{K}(x,y) &\leq (x-y)^{1-1/p}\Big(\int_y^x t^{p-1}|v'(t)|^p dt\Big)^{1/p}\\
\text{(by (2.4))} &\leq c_1(p)(x-y)^{1-1/p}(I[u])^{1/p}\\
&\leq c_1(p)(x-y)^{1-1/p}X^{-1/p}((x-y)/D)(I[u])^{1/p},
\end{aligned} \qquad (6.3)$$

for any $D \geq R$, since $0 \leq X(t) \leq 1$ for all $t \in (0,1]$. To estimate $\Lambda(x,y)$ we return to the original variable by $v(t) = u(t)/t^{1-1/p}$, thus

$$\Lambda(x,y) = \int_y^x \frac{|u(t)|}{t} dt.$$

Inserting (6.1) in $\Lambda(x,y)$ we obtain

$$\begin{aligned}
\Lambda(x,y) &\leq c_2(p)(I[u])^{1/p}\int_y^x t^{-1/p}X^{-1/p}(t/D)dt\\
&\leq c_3(p)(x-y)^{1-1/p}X^{-1/p}((x-y)/D)(I[u])^{1/p},
\end{aligned} \qquad (6.4)$$

for any $D \geq e^\eta R$, where $\eta$ depends only on $p$, due to Lemma 2.5 $(ii)$ for $\alpha = -1/p$ and $\beta = 1/p$. Coupling (6.3) and (6.4) with (6.2), we end up with

$$|u(x) - u(y)| \leq c_4(p)(x-y)^{1-1/p}X^{-1/p}((x-y)/D)(I[u])^{1/p},$$

for all $0 < y < x < R$ and any $D \geq e^\eta R$, which is the desired estimate with $B = e^\eta$. ∎

# 7 Optimality of the logarithmic correction

In this section we prove the optimality of the exponent $1/p$ on $X$, in the Hölder semi-norm inequality (1.8) of Theorem B. Note that we can pick one point in (1.8) to be the origin, and therefore it is enough to prove the alleged optimality in (5.1) and (6.1).

We consider the radially symmetric, Lipschitz continuous function

$$u_\delta(x) = \begin{cases}
(\delta^2|x|)^H(6 - \frac{\log|x|}{\log\delta}), & \delta^6 \leq |x| < \delta^5\\
(\delta^{-3}|x|^2)^H, & \delta^5 \leq |x| < \delta^4\\
(\delta|x|)^H(1 + H\log(|x|/\delta^4)), & \delta^4 \leq |x| < \delta^3\\
(\delta|x|)^H(1 - H\log(|x|/\delta^2)), & \delta^3 \leq |x| < \delta^2\\
\delta^{3H}, & \delta^2 \leq |x| < \delta\\
(\delta^2|x|)^H \frac{\log|x|}{\log\delta}, & \delta \leq |x| \leq 1
\end{cases}$$

where $0 < \delta < 1$ and $H := (p-n)/p$ with $p > n \geq 1$. With $u_\delta$ we associate the quotient

$$Q_\epsilon[u_\delta;x] := \frac{(I[u_\delta])^{1/p}}{|u_\delta(x)||x|^{-H}X^{1/p-\epsilon}(|x|/D)}, \quad 0 \leq \epsilon < 1/p, \quad \delta^6 < |x| < 1.$$

Note that due to (5.1) and (6.1) (and after an approximation of $u_\delta$ by smooth functions), we have $Q_0[u_\delta; x] \geq C$, for some positive constant $C = C(n, p)$. To prove that the exponent $1/p$ on the correction weight $X$ cannot be decreased, we fix $0 < \epsilon < 1/p$ and taking $x_\delta$ such that $|x_\delta| = \delta^3$ we will prove that $Q_\epsilon[u_\delta; x_\delta] \to 0$ as $\delta \downarrow 0$.

We begin by computing $I[u_\delta]$. Setting $A_k := \{x \in \mathbb{R}^n : \delta^k < |x| < \delta^{k-1}\}$, $k = 1, ..., 6$, we have $I[u_\delta] = \sum_{k=1}^{6} I[u_\delta](A_k)$, where

$$I[u_\delta](A_k) := \int_{A_k} |\nabla u_\delta(x)|^p dx - H^p \int_{A_k} \frac{|u_\delta(x)|^p}{|x|^p} dx, \quad k = 1, ..., 6.$$

By the fact that $u_\delta$ is radially symmetric, we may use polar coordinates to get

$$I[u_\delta](A_k) = n\omega_n \left[ \int_{\delta^k}^{\delta^{k-1}} |\tilde{u}'_\delta(t)|^p t^{n-1} dt - H^p \int_{\delta^k}^{\delta^{k-1}} |\tilde{u}_\delta(t)|^p t^{n-1-p} dt \right], \quad k = 1, ..., 6,$$

where $\tilde{u}_\delta(t) = u_\delta(x)$ with $t = |x|$. We then have

$$I[u_\delta](A_1) = n\omega_n \frac{\delta^{2pH}}{\log^p(1/\delta)} \left[ \int_{\delta}^{1} t^{-1} |1 - H\log(1/t)|^p dt - \int_{\delta}^{1} t^{-1} (H\log(1/t))^p dt \right]$$

$$= n\omega_n \frac{\delta^{2pH}}{\log^p(1/\delta)} \left[ \int_{\delta}^{e^{-1/H}} + \int_{e^{-1/H}}^{1} t^{-1} |1 - H\log(1/t)|^p dt - \int_{\delta}^{1} t^{-1} (H\log(1/t))^p dt \right]$$

$$= \frac{n\omega_n}{(p+1)H} \delta^{2pH} \log(1/\delta) \left[ \left( H - \frac{1}{\log(1/\delta)} \right)^{p+1} + \left( \frac{1}{\log(1/\delta)} \right)^{p+1} - H^{p+1} \right],$$

where (since we will let $\delta \downarrow 0$) we have assumed $\delta < e^{-1/H}$ in order to get rid of the absolute value. Now we compute

$$I[u_\delta](A_6) = n\omega_n \frac{\delta^{2pH}}{\log^p(1/\delta)} \left[ \int_{\delta^6}^{\delta^5} t^{-1} (1 + H\log(t/\delta^6))^p dt - \int_{\delta^6}^{\delta^5} t^{-1} (H\log(t/\delta^6))^p dt \right]$$

$$= \frac{n\omega_n}{(p+1)H} \delta^{2pH} \log(1/\delta) \left[ \left( H + \frac{1}{\log(1/\delta)} \right)^{p+1} - \left( \frac{1}{\log(1/\delta)} \right)^{p+1} - H^{p+1} \right].$$

Thus $I[u_\delta](A_1) + I[u_\delta](A_6) =$

$$\frac{2n\omega_n}{(p+1)H} \delta^{2pH} \log(1/\delta) \left[ \left( H + \frac{1}{\log(1/\delta)} \right)^{p+1} + \left( H - \frac{1}{\log(1/\delta)} \right)^{p+1} - 2H^{p+1} \right].$$

The factor in the square brackets is of order $o(1)$, as $\delta \downarrow 0$. Since $H = (p-n)/p$, we get

$$I[u_\delta](A_1) + I[u_\delta](A_6) = o(\delta^{2(p-n)} \log(1/\delta)), \quad \text{as } \delta \downarrow 0. \tag{7.1}$$

Similarly,

$$I[u_\delta](A_2) = -n\omega_n H^p \delta^{3pH} \int_{\delta^2}^{\delta} t^{-pH-1} dt$$

$$= -\frac{n\omega_n}{p} H^{p-1} \delta^{pH} (1 - \delta^{pH}),$$

and

$$I[u_\delta](A_5) = n\omega_n H^p \delta^{-3pH}\left[2^p \int_{\delta^5}^{\delta^4} t^{pH-1}dt - \int_{\delta^5}^{\delta^4} t^{pH-1}dt\right]$$

$$= \frac{n\omega_n}{p}H^{p-1}(2^p - 1)\delta^{pH}(1 - \delta^{pH}).$$

Hence

$$I[u_\delta](A_2) + I[u_\delta](A_5) = \frac{n\omega_n}{p}H^{p-1}(2^p - 2)\delta^{pH}(1 - \delta^{pH})$$
$$= O(\delta^{p-n}), \quad \text{as } \delta \downarrow 0. \tag{7.2}$$

Finally, the first summand of the last pair is

$$I[u_\delta](A_3) = n\omega_n H^p \delta^{pH}\left[\int_{\delta^3}^{\delta^2}\left(-H\log\frac{t}{\delta^2}\right)^p t^{-1}dt - \int_{\delta^3}^{\delta^2}\left(1 - H\log\frac{t}{\delta^2}\right)^p t^{-1}dt\right]$$

$$= \frac{n\omega_n}{p+1}H^{2p}\delta^{pH}\left(\log\frac{1}{\delta}\right)^{p+1}\left[1 + \frac{1}{(H\log\frac{1}{\delta})^{p+1}} - \left(1 + \frac{1}{H\log\frac{1}{\delta}}\right)^{p+1}\right],$$

and the second one

$$I[u_\delta](A_4) = n\omega_n H^p \delta^{pH}\left[\int_{\delta^4}^{\delta^3}\left(2 + H\log\frac{t}{\delta^4}\right)^p t^{-1}dt - \int_{\delta^4}^{\delta^3}\left(1 + H\log\frac{t}{\delta^4}\right)^p t^{-1}dt\right]$$

$$= \frac{n\omega_n}{p+1}H^{2p}\delta^{pH}\left(\log\frac{1}{\delta}\right)^{p+1}\left[\left(1 + \frac{2}{H\log\frac{1}{\delta}}\right)^{p+1} - \left(1 + \frac{1}{H\log\frac{1}{\delta}}\right)^{p+1} + \frac{1 - 2^{p+1}}{(H\log\frac{1}{\delta})^{p+1}}\right].$$

Adding, we find $I[u_\delta](A_3) + I[u_\delta](A_4) =$

$$\frac{n\omega_n}{p+1}H^{2p}\delta^{pH}\left(\log\frac{1}{\delta}\right)^{p+1}\left[1 + \left(1 + \frac{2}{H\log\frac{1}{\delta}}\right)^{p+1} - 2\left(1 + \frac{1}{H\log\frac{1}{\delta}}\right)^{p+1} - \frac{2^{p+1} - 2}{(H\log\frac{1}{\delta})^{p+1}}\right].$$

The factor in square brackets is of order $\frac{p(p+1)}{H^2}o(\frac{1}{\log^2(1/\delta)})$, as $\delta \downarrow 0$, and we get

$$I[u_\delta](A_3) + I[u_\delta](A_4) = pn\omega_n H^{2(p-1)}\delta^{p-n}(\log(1/\delta))^{p-1} + o(\delta^{p-n}(\log(1/\delta))^{p-1}), \tag{7.3}$$

as $\delta \downarrow 0$. From (7.1), (7.2) and (7.3), we see that the leading term in $I[u_\delta] = \sum_{k=1}^{6} I[u_\delta](A_k)$, comes from $I[u_\delta](A_3) + I[u_\delta](A_4)$. More precisely

$$I[u_\delta] = pn\omega_n H^{2(p-1)}\delta^{p-n}(\log(1/\delta))^{p-1} + o(\delta^{p-n}(\log(1/\delta))^{p-1}),$$

as $\delta \downarrow 0$. Finally, we compute the denominator of $Q_\epsilon[u_\delta; x_\delta]$ with $|x_\delta| = \delta^3$,

$$|u_\delta(\delta^3)|\delta^{-3H}X^{1/p-\epsilon}(\delta^3/D) = \delta^H(1 + H\log(1/\delta))X^{1/p-\epsilon}(\delta^3/D)$$
$$= \delta^H(1 + H\log(1/\delta))(1 - \log(\delta^3/D))^{-1/p+\epsilon}$$
$$= 3H\delta^H(\log(1/\delta))^{1-1/p+\epsilon} + o(\delta^H(\log(1/\delta))^{1-1/p}),$$

as $\delta \downarrow 0$. Dividing the two endmost relations we conclude

$$\begin{aligned}
Q_\epsilon[u_\delta; x_\delta] &= \frac{1}{3}(pn\omega_n)^{1/p}H^{1-2/p}\frac{[\delta^{p-n}(\log 1/\delta)^{p-1} + o(\delta^{p-n}(\log(1/\delta))^{p-1})]^{1/p}}{\delta^H(\log(1/\delta))^{1-1/p+\epsilon} + o(\delta^{1-n/p}(\log(1/\delta))^{1-1/p})} \\
&= \frac{1}{3}(pn\omega_n)^{1/p}H^{1-2/p}\frac{[1+o(1)]^{1/p}}{(\log(1/\delta))^\epsilon + o(1)} \\
&\to 0, \quad \text{as } \delta \downarrow 0.
\end{aligned}$$
∎

## 8 Optimality of the logarithmic correction II

In this section we give another proof of the optimality of the exponent $1/p$ on $X$, in the Hardy-Morrey inequality of Theorem B.

Recall first the following technical lemma from [BFT] (Lemma 5.2 (i)).

**Lemma 8.1.** *For $\beta > -1$, small $\varepsilon > 0$, $p > 1$ and $\delta \leq 1$, set $J_\beta(\varepsilon) := \int_0^\delta t^{-1+\varepsilon p} X^{-\beta}(t)dt$. Then $J_\beta(\varepsilon) \asymp \varepsilon^{-1-\beta}$.*

To prove that the exponent $1/p$ on $X$ in (1.8) is optimal, for arbitrarily small $\varepsilon > 0$ and $0 < \theta < 1/p$, we introduce the function

$$U_\varepsilon(x) := |x|^{H+\varepsilon} X^{-\theta}(|x|); \quad x \in B_1(0) \setminus \{0\}.$$

Let also $B_\delta$ be a ball of radius $0 < \delta < \min\{e^{1-\theta/H}, 1\}$, centered at the origin and such that $B_\delta \subset\subset \Omega$. We consider a radially symmetric $\phi \in C_c^\infty(B_\delta)$ with the properties: (i) $0 \leq \phi(x) \leq 1$; $x \in B_\delta$, (ii) $\phi(x) \equiv 1$; $x \in B_{\delta/2}$, and (iii) $|\nabla\phi(x)| \leq 1/\delta$; $x \in B_\delta$. Now we define

$$u_\varepsilon(x) := U_\varepsilon(x)\phi(x); \quad x \in \Omega.$$

Obviously, $\text{sprt}\{u_\varepsilon\} \subset B_\delta$ and $\lim_{x\to 0} u_\varepsilon(x) = 0$. Also

$$\nabla u_\varepsilon = \phi \nabla U_\varepsilon + U_\varepsilon \nabla \phi(x),$$

and using the elementary inequality

$$|\alpha + \beta|^p \leq |\alpha|^p + c_p\Big(|\alpha|^{p-1}|\beta| + |\beta|^p\Big); \quad \alpha, \beta \in \mathbb{R}^n,$$

we obtain

$$\int_\Omega |\nabla u_\varepsilon(x)|^p dx \leq$$
$$\underbrace{\int_{B_\delta} |\phi(x)\nabla U_\varepsilon(x)|^p dx}_{=:J_1} + c_p \underbrace{\int_{B_\delta} |\phi(x)\nabla U_\varepsilon(x)|^{p-1}|U_\varepsilon \nabla \phi(x)|dx}_{=:J_2} + c_p \underbrace{\int_{B_\delta} |U_\varepsilon \nabla \phi(x)|^p dx}_{=:J_3}.$$

Simple calculations and taking into account the properties we imposed on $\phi$, give

$$\begin{aligned} J_3 &\leq \delta^{-p} \int_{B_\delta} |U_\varepsilon|^p dx \\ &= n\omega_n \delta^{-p} \int_0^\delta t^{p-1+\varepsilon p} X^{-\theta p}(t) dt \\ &= O_\varepsilon(1), \quad \text{as } \varepsilon \downarrow 0, \end{aligned}$$

while

$$\begin{aligned} J_2 &\leq \delta^{-1} \int_{B_\delta} |\nabla U_\varepsilon(x)|^{p-1} |U_\varepsilon| dx \\ &= n\omega_n \delta^{-1} \int_0^\delta t^{\varepsilon p} X^{-\theta p}(t) |H + \varepsilon - \theta X(t)|^{p-1} dt \\ &\leq n\omega_n \delta^{-1} (H+\varepsilon)^{p-1} \int_0^\delta t^{\varepsilon p} X^{-\theta p}(t) dt \\ &= O_\varepsilon(1), \quad \text{as } \varepsilon \downarrow 0. \end{aligned}$$

Here we have used the fact that since $\delta < e^{1-\theta/H}$, we have $H + \varepsilon - \theta X(t) \geq 0$, for all $t \in [0,\delta]$ and all $\varepsilon \geq 0$. Next, by the above considerations we obtain

$$\begin{aligned} I[u_\varepsilon] &\leq J_1 - H^p \int_\Omega \frac{|u_\varepsilon(x)|^p}{|x|^p} dx + O_\varepsilon(1) \\ &= \int_{B_\delta} \phi^p(x) \Big[ |\nabla U_\varepsilon(x)|^p - H^p \frac{U_\varepsilon^p(x)}{|x|^p} \Big] dx + O_\varepsilon(1) \\ &= n\omega_n \int_0^\delta \phi^p(t) t^{-1+\varepsilon p} X^{-\theta p}(t) \Big\{ \Big(H + \varepsilon - \theta X(t)\Big)^p - H^p \Big\} dt + O_\varepsilon(1), \end{aligned}$$

as $\varepsilon \downarrow 0$. For small enough $\varepsilon, \theta$ we have $|\varepsilon - \theta X(t)| << H$ so that we can write

$$\Big(H+\varepsilon-\theta X(t)\Big)^p - H^p \leq pH^{p-1}\Big(\varepsilon-\theta X(t)\Big) + \frac{1}{2}p(p-1)H^{p-2}\Big(\varepsilon-\theta X(t)\Big)^2 + c\Big(\varepsilon-\theta X(t)\Big)^3,$$

by Taylor's expansion theorem. Thus, as $\varepsilon \downarrow 0$,

$$\begin{aligned} I[u_\varepsilon] &\leq n\omega_n p H^{p-1} \int_0^\delta \phi^p(t) t^{-1+\varepsilon p} X^{-\theta p}(t) \Big(\varepsilon - \theta X(t)\Big) dt \\ &\quad + \frac{1}{2} n\omega_n p(p-1) H^{p-2} \int_0^\delta \phi^p(t) t^{-1+\varepsilon p} X^{-\theta p}(t) \Big(\varepsilon - \theta X(t)\Big)^2 dt \\ &\quad + n\omega_n c \int_0^\delta \phi^p(t) t^{-1+\varepsilon p} X^{-\theta p}(t) \Big(\varepsilon - \theta X(t)\Big)^3 dt + O_\varepsilon(1), \\ &=: J_{11} + J_{12} + J_{13} + O_\varepsilon(1). \end{aligned}$$

Noting that

$$\varepsilon \int_0^\delta \phi^p(t) t^{-1+\varepsilon p} X^{-\theta p}(t) dt = \frac{1}{p} \int_0^\delta \phi^p(t) (t^{\varepsilon p})' X^{-\theta p}(t) dt$$
$$= \theta \int_0^\delta \phi^p(t) t^{-1+\varepsilon p} X^{1-\theta p}(t) dt$$
$$- \int_0^\delta \phi^{p-1}(t) \phi'(t) t^{\varepsilon p} X^{-\theta p}(t) dt,$$

we deduce

$$J_{11} = -n\omega_n p H^{p-1} \int_0^\delta \phi^{p-1}(t) \phi'(t) t^{\varepsilon p} X^{-\theta p}(t) dt$$
$$\leq n\omega_n p H^{p-1} \delta^{-1} \int_0^\delta t^{\varepsilon p} X^{-\theta p}(t) dt$$
$$= O_\varepsilon(1).$$

Using Lemma 8.1, it is easy to verify that $J_{12}, J_{13} = o_\varepsilon(1)$, as $\varepsilon \downarrow 0$. For example, for $J_{12}$ we have $\varepsilon^2 J_{\theta p}(\varepsilon), \varepsilon J_{\theta p-1}(\varepsilon), J_{\theta p-2}(\varepsilon) \asymp \varepsilon^{1-\theta p} \to 0$, since $0 < \theta < 1/p$. We conclude that $I[u_\varepsilon] = O_\varepsilon(1)$, as $\varepsilon \downarrow 0$.

Now assume that for some $\epsilon \in (0, 1/p]$ we have

$$\sup_{\substack{x,y \in \Omega \\ x \neq y}} \left\{ \frac{|u(x) - u(y)|}{|x-y|^{1-n/p}} X^{1/p-\epsilon}\left(\frac{|x-y|}{D}\right) \right\} \leq C(I[u])^{1/p},$$

for all $u \in W_0^{1,p}(\Omega \setminus \{0\})$. Then for $u = u_\varepsilon$ with $\theta = 1/p - \epsilon/2$, and $x \in B_{\delta/2}(0), y = 0$, we obtain

$$|x|^\varepsilon X^{-\epsilon/2}(|x|/D) \leq C(I[u_\varepsilon])^{1/p}.$$

Letting $\varepsilon \downarrow 0$, we obtain $X^{-\epsilon/2}(|x|/D) \leq C$, for all $x \in B_{\delta/2}(0)$, which is absurd. ∎

**Acknowledgements** The author would like to thank Prof. Stathis Filippas for his essential help in this work.